\documentclass[notitlepage,11pt,a4paper]{scrartcl}
\usepackage{amstext}
\usepackage{amssymb}
\usepackage{amsfonts}
\usepackage{amssymb}
\usepackage{amsmath}
\usepackage{graphicx}
\usepackage[all]{xy}
\usepackage[latin1]{inputenc}
\usepackage{theorem}
\usepackage{geometry}
\usepackage{abstract}
\newtheorem{satz}{Theorem}[section]

\newtheorem{defi}[satz]{Definition}
\newtheorem{bem}[satz]{Remark}
\newtheorem{kor}[satz]{Corollary}
\newtheorem{lem}[satz]{Lemma}

\newtheorem{bei}[satz]{Example}
\newtheorem{pro}[satz]{Proposition}

\newcommand{\qed}{\begin{flushright}$\square$\end{flushright}}
\newcommand{\filt}[6]{
\[
\begin{xy}
\xymatrix@R20pt@C20pt{
&\mathbb{C}^3&\\\langle #1,#2 \rangle\ar[ru]&\langle #3,#4\rangle\ar[u]&\langle #5,#6\rangle\ar[lu]\\\langle #1\rangle\ar[u]&\langle #3\rangle\ar[u]&\langle #5\rangle\ar[u]}
\end{xy}
\]
}

\newcommand\Zn{\mathbb{Z}}
\newcommand\Nn{\mathbb{N}}

\newcommand{\Sc}[2]{\langle #1,#2\rangle}
\newcommand{\Hom}{\mathrm{Hom}} 
\newcommand{\Ext}{\mathrm{Ext}} 
\newcommand{\ext}{\mathrm{ext}} 
\newcommand{\End}{\mathrm{End}} 
\newcommand{\Rep}{\mathrm{Rep}} 
\newcommand{\ses}[3]{0\rightarrow #1\rightarrow #2\rightarrow#3\rightarrow 0}
\begin{document}
  \title{Tree modules}
  \author{Thorsten Weist\\Fachbereich C - Mathematik\\
Bergische Universität Wuppertal\\
D - 42097 Wuppertal, Germany\\
\texttt{e-mail: weist@math.uni-wuppertal.de}}
\date{}   
\maketitle
\begin{abstract}
\noindent After stating several tools which can be used to construct indecomposable tree modules for quivers without oriented cycles, we use these methods to construct indecomposable tree modules for every imaginary Schur root. These methods also give a recipe for the construction of tree modules for every root. Moreover, we give several examples illustrating the results.
\end{abstract}
%{\bf\LARGE Tree modules
%\\}
%\end{center}
%\vspace{0.0cm}
%\begin{center}
%Thorsten Weist
%\end{center}
%\vspace{0.0cm}

\section{Introduction}
In this paper we study tree modules of hereditary finite dimensional $k$-algebras or equivalently of quivers without relations and oriented cycles. Therefore, given a representation of a quiver, we choose a basis of the vector spaces associated to each vertex of the quiver and consider the maps restricted to these basis elements. Then we investigate the {\it coefficient quiver} in which the basis elements label the vertices and which has an arrow between two vertices if the matrix coefficient corresponding to these two basis elements is not zero. A representation is called a {\it tree module} if there exists a basis such that the coefficient quiver is a tree. This leads us to the following problem stated by Ringel, see \cite[Problem 9]{rin3}:\\
Does there exist an indecomposable tree module for every wild hereditary quiver and every root? In particular, Ringel conjectured that there should be more than one isomorphism class for imaginary roots.\\
The main result of this paper is that there exists more than one isomorphism class of indecomposable tree modules for every isotropic root and for every imaginary Schur root, see Theorems \ref{isot} and \ref{imag}. In the course of the proof we state explicit methods which describe how to construct these indecomposable tree modules. Moreover, the construction is also applicable for many non-Schurian roots. These methods can also be used to construct indecomposable modules which are not necessarily tree modules.\\
As far as Schur roots are concerned, we determine an exceptional sequence of real Schur roots corresponding to a fixed imaginary root. By \cite{rin1} the unique indecomposable representations of these roots are known to be trees. These representations are the building blocks of certain indecomposable representations of the imaginary Schur root. If the sequence has length two, all representations which can be built by these two correspond to representations of the (generalised) Kronecker quiver, i.e. the quiver with two vertices and no oriented cycles. For this quiver the conjecture is known to be true, see \cite{wei}. If the sequence has length greater than two, we can either recursively apply Ringel's reflection functors or glue appropriate tree modules of smaller dimension in order to construct indecomposable tree modules for every Schur root.\\
In the last section we consider several examples illustrating the results of the preceding sections. Therefore, it is often useful to consider the universal cover of the given quiver. This quiver is already a tree and, moreover, it is known that the functor which maps the representations of the universal cover to the original one preserves indecomposability. Thus even if the methods of the paper are applicable without restrictions to the shape of the quiver (except those we made in the beginning), for explicit calculations it is usually more convenient to consider the universal covering quiver. Moreover, every tree module is already a representation of the universal covering quiver.\\
In order to decide whether a given dimension vector is a Schur root or not, one can consider the canonical decomposition of the dimension vector. But not much is known about the ratio of Schurian to non-Schurian roots in general.  One example is discussed in \cite{dw3} where the imaginary Schur roots of the considered quiver are given by a non-convex fractal-like polygon contained in a quadric describing all imaginary roots. Moreover, based on the algorithm of \cite{dw} by the methods presented in this paper it is possible to construct Schurian representations and, therefore, also Schur roots. In conclusion, it is hard to say how much of Ringel's conjecture is still an open problem at this point, also because the stated methods are applicable for many, but not all, non-Schurian roots as well. We should also note that the known examples for which the methods do not apply, see \cite{wie} and Example \ref{bei}, give the impression that these roots are constructable in a certain way, but that the majority of non-Schur roots is not of this type.

\section{Generalities}\label{allg}
Let $k$ be an algebraically closed field.
\begin{defi}
A quiver $Q$ consists of a set of vertices $Q_0$ and a set of arrows $Q_1$ denoted by $\rho:i\rightarrow j$ for $i,j\in Q_0$.\\A vertex $i\in Q_0$ is called sink if there does not exist an arrow $\rho:i\rightarrow j\in Q_1$.\\
A vertex $j\in Q_0$ is called source if there does not exist an arrow $\rho:i\rightarrow j\in Q_1$.
\end{defi}
Define the abelian group
\[\mathbb{Z}Q_0=\bigoplus_{i\in Q_0}\mathbb{Z}i\] and its monoid of dimension vectors $\mathbb{N}Q_0$.\\
A finite-dimensional $k$-representation of $Q$ is given by a tuple
\[X=((X_i)_{i\in Q_0},(X_{\rho})_{\rho\in Q_1}:X_i\rightarrow X_j)\]
of finite-dimensional $k$-vector spaces and $k$-linear maps between them. The dimension vector $\underline{\dim}X\in\mathbb{N}Q_0$ of $X$ is defined by
\[\underline{\dim}X=\sum_{i\in Q_0}\dim_kX_ii.\]
A dimension vector is called a root if there exists an indecomposable representation of this dimension vector. It is called Schur root if there exists a representation with trivial endomorphism ring.\\
Let $\alpha\in\mathbb{N}Q_0$ be a dimension vector. The variety $R_\alpha(Q)$ of $k$-representations of $Q$ of dimension vector
$\alpha$ is defined as the affine $k$-space
\[R_{\alpha}(Q)=\bigoplus_{\rho:i\rightarrow j} \mathrm{Hom}_k(k^{\alpha_i},k^{\alpha_j}).\]
%The algebraic group 
%\[G_d=\prod_{i\in I} Gl_{d_i}(k)\]
%acts on $R_d(Q)$ via simultaneous base change, i.e.
%\[(g_i)_{i\in Q_0}\ast
%(X_{\alpha})_{\alpha\in Q_1}=(g_jX_{\alpha}g_i^{-1})_{\alpha:i\rightarrow
%j}.\] The orbits are in bijection with the isomorphism classes of
%$k$-representations of $Q$ with dimension vector $d$.\\\\ 
In the space of $\mathbb{Z}$-linear functions $\mathrm{Hom}_{\mathbb{Z}}(\mathbb{Z}Q_0,\mathbb{Z})$ we consider the basis given by the elements $i^{\ast}$ for $i\in Q_0$, i.e.
$i^{\ast}(j)=\delta_{i,j}$ for $j\in Q_0$. Define
\[\dim:=\sum_{i\in Q_0}i^{\ast}.\]
On $\Zn Q_0$ we have a non-symmetric bilinear form, the Euler form,
which is defined by
\[\Sc{\alpha}{\beta}=\sum_{i\in Q_0}\alpha_i\beta_i-\sum_{\rho:i\rightarrow j\in Q_1}\alpha_i\beta_j\]
for $\alpha,\,\beta\in\Zn Q_0$.\\
By \cite{rin2}, for two representations $X$, $Y$ of $Q$ we have
\[\Sc{\underline{\dim} X}{\underline{\dim} Y}=\dim_k\Hom(X,Y)-\dim_k\Ext(X,Y)\]
and $\Ext^i_k(X,Y)=0$ for $i\geq 2$.\\
If some property is independent of the point chosen in some open subset $U$ of $R_{\alpha}(Q)$, following \cite{sch}, we say that this property is true for a general representation of dimension vector $\alpha\in\Nn Q_0$.\\
\\Since the function
$\lambda:R_{\alpha}(Q)\times R_{\beta}(Q)\rightarrow\mathbb{N},\,(X,Y)\mapsto\dim\Hom(X,Y)$ is upper semi-continuous, see for instance \cite{sch}, we can define $\hom(\alpha,\beta)$ to be the minimal, and therefore general, value of this function. In particular, we get that if $\alpha$ is a Schur root of a quiver, then a general representation is Schurian. Moreover, let $\ext(\alpha,\beta):=\hom(\alpha,\beta)-\Sc{\alpha}{\beta}$.\\
We denote by $\beta\hookrightarrow \alpha$ if a general representation of dimension $\alpha$ has a subrepresentation of dimension $\beta$.\\\\
Following \cite{kac}, for every dimension vector $\alpha\in\Nn Q_0$ there exists a decomposition $\alpha=\sum_{i\in I}\beta_i$ and an open subset of $R_{\alpha}(Q)$ such that a general representation $Y\in U$ is a direct sum of Schurian representation $X_i$ with $\underline{\dim} X_i=\beta_i$. We write $\alpha=\oplus_{i\in I}\beta_i$. This is called the canonical decomposition of $\alpha$. Moreover, we have the following result, see \cite{kac} and \cite[Theorem 4.4]{sch}:
\begin{satz}
\begin{enumerate}
\item For a general representation $Y$ of dimension vector $\alpha$ we have $Y\cong\oplus_{i\in I} X_i$ with $\underline{\dim} X_i=\beta_i$ if and only if $\ext(\beta_i,\beta_j)=0$ for $i\neq j$. Moreover, each $X_i$ is Schurian.
\item Let $\alpha$ be a root. Then up to multiplicity there exists at most one imaginary Schur root in its decomposition.
\end{enumerate}
\end{satz}
Note that \cite[Section 4]{dw} gives a very useful algorithm which can be used to determine the canonical decomposition.\\ 
\setcounter{subsection}{1}
We introduce coefficient quivers and tree modules following the presentation given in $\cite{rin1}$.\\
Let $Q$ be a quiver, $\alpha=(\alpha_q)_{q\in Q_0}$ a dimension vector and $X$ with $\underline{\dim} X=\alpha$ a representation of $Q$. A basis of $X$ is a subset $\mathcal{B}$ of $\bigoplus_{q\in Q_0}X_q$ such that
\[\mathcal{B}_q:=\mathcal{B}\cap X_q\] is a basis of $X_q$ for all vertices $q\in Q_0$. For every arrow $\rho:i\rightarrow j$ we may write $X_{\rho}$ as a $(\alpha_j\times \alpha_i)$-matrix $X_{\rho,\mathcal{B}}$ with coefficients in $k$ such that the rows and columns are indexed by $\mathcal{B}_j$ and $\mathcal{B}_i$ respectively. If
\[X_{\rho}(b)=\sum_{b'\in\mathcal{B}_j}\lambda_{b',b}b'\]
with $\lambda_{b',b}\in k$, we obviously have $(X_{\rho,\mathcal{B}})_{b',b}=\lambda_{b',b}$.
\begin{defi}
The coefficient quiver $\Gamma(X,\mathcal{B})$ of a representation $X$ with a fixed basis $\mathcal{B}$ has vertex set $\mathcal{B}$ and arrows between vertices are defined by the condition: if $(X_{\rho,\mathcal{B}})_{b,b'}\neq 0$, there exists an arrow $(\rho,b,b'):b\rightarrow b'$.\\
A representation $X$ is called a tree module if there exists a basis $\mathcal{B}$ for $X$ such that the corresponding coefficient quiver is a tree.
\end{defi}
%As already mentioned in the introduction, defining the coefficient quiver immediately raises the %question posed by Claus Michael Ringel:\\
%Does there exist an indecomposable tree module for every wild hereditary quiver and every root $d$? He conjectured that there should be more than one isomorphism class for imaginary roots.\\\\
In order to construct a tree module and its coefficient quiver respectively, it is often useful to consider the universal covering quiver of the given quiver $Q=(Q_0,Q_1)$.\\
Let $Q=(Q_0,Q_1)$ be a connected quiver without oriented cycles. Let $Q^{-1}_1=\{\rho,\rho^{-1}\mid\rho\in Q_1\}$ where $\rho^{-1}$ is the formal inverse of $\rho$. We will write $\rho^{-1}:j\rightarrow i$ for $\rho:i\rightarrow j\in Q_1$. A path $p$ is a sequence $(i_1\mid\rho_1\rho_2\ldots\rho_n\mid i_{n+1})$ such that $\rho_j:i_j\rightarrow i_{j+1}\in Q^{-1}_1$. Thereby, we have the equivalence generated by
\[(i\mid\rho\rho^{-1}\mid i)\sim (i\mid\mid i).\]
In what follows, we always consider paths up to this equivalence. The set of words in $Q$ is generated by the arrows and their formal inverses, i.e. for a word $w$ we have $w=\rho_1\ldots\rho_n$ where $\rho_i\in Q_1^{-1}$. Denote the set of words of $Q$ by $W(Q)$. The universal covering quiver $\tilde{Q}$ of $Q$ is given by the vertex set 
\[\tilde{Q}_0=\{(i,w)\mid i\in Q_0,w\in W(Q)\}\]
and the arrow set
\[\tilde{Q}_1=\{\rho_{(i,w)}:(i,w)\rightarrow (j,w\rho)\mid\rho:i\rightarrow j\in Q_1\}.\]
Every representation $\tilde{X}$ of $\tilde{Q}$ gives rise to a representations of $X$ of $Q$ in the following way:
\[X_i=\bigoplus_{w\in W(Q)} \tilde{X}_{(i,w)},\quad i\in Q_0\]
and $X_{\rho}:X_i\rightarrow X_j$ is defined by $X_{\rho}|_{X_{(i,w)}}=\tilde{X}_{\rho_{(i,w)}}$. Now we can make use of the following result, see \cite[Lemma 3.5]{gab} or \cite[Theorem 2.2]{gre}:
\begin{satz}
If $\tilde{X}$ is an indecomposable representation of $\tilde{Q}$, the corresponding representation $X$ of $Q$ is also indecomposable.
\end{satz}
Note, that the endomorphism rings of $X$ and $\tilde{X}$ do not have to coincide. But we clearly have $\End(\tilde{X})\subseteq\End(X)$.\\\\
It is straightforward to check that every indecomposable tree module is already a representation of a connected component of the universal cover. Indeed, start by fixing some vertex $i\in\mathcal{B}$ of the coefficient quiver. Since the tree is connected, every vertex $j\in\mathcal{B}$ defines a path from $i$ to $j$. Moreover, every vertex of the coefficient quiver corresponds to a vertex of the original quiver.
\begin{bei}
\end{bei}
Consider the quiver
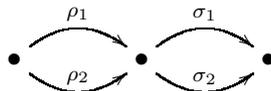
\[
\begin{xy}
\xymatrix@R20pt@C15pt{\bullet\ar@/^1.0pc/[rr]^{\rho_1}\ar@/_1.0pc/[rr]^{\rho_2}&&\bullet\ar@/^1.0pc/[rr]^{\sigma_1}\ar@/_1.0pc/[rr]^{\sigma_2}&&\bullet
}
\end{xy}\] 
and the real root $(1,2,4)$. Now it is easy to write down the indecomposable tree module of the universal cover corresponding to this dimension vector, i.e.
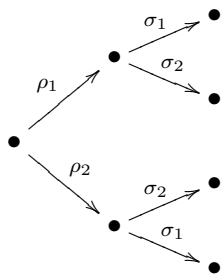
\[
\begin{xy}
\xymatrix@R5pt@C10pt{&&&&\bullet\\&&\bullet\ar[rru]^{\sigma_1}\ar[rrd]^{\sigma_2}\\&&&&\bullet\\\bullet\ar[rruu]^{\rho_1}\ar[rrdd]^{\rho_2}\\&&&&\bullet\\
&&\bullet\ar[rrd]^{\sigma_1}\ar[rru]^{\sigma_2}\\&&&&\bullet
}
\end{xy}\] 
Here the dots represent vector spaces of dimension one and the arrows the identity map.

\section{Basic tools}
In this section we present different tools which can be used to construct tree modules. Therefore, we first consider the reflection functor introduced in \cite{rin} and combine it with several results of \cite{sc2} dealing with perpendicular categories. Afterwards we recall some results of \cite{sch} concerning the canonical decomposition of the dimension vectors of a quiver. Moreover, we need the main result of \cite{wei}, i.e. the existence of indecomposable tree modules for every root of the Kronecker quiver.
\subsection{Exceptional sequences and Reflection functors}
We denote by $\mathrm{Rep}(Q)$ the category of finite-dimensional representations of $Q$. An indecomposable representation $X$ of a quiver $Q$ is called exceptional if $\Ext(X,X)=0$. Then it follows that $\underline{\dim} X$ is a real Schur root and $\mathrm{End}(X)=k$, see also Lemma \ref{happelringel}.\\ A sequence $S=(X_1,\ldots,X_r)$ of representations of $Q$ is called exceptional if every $X_i$ is exceptional and, moreover, $\Hom(X_i,X_j)=\Ext(X_i,X_j)=0$ if $i<j$. If we do not require that the representations $X_i$ are exceptional, such a sequence is called Schur sequence.\\
For a set $S=\{X_1,\ldots,X_r\}$ of representations of $Q$ we define its perpendicular categories
\[^{\perp}S=\{X\in\mathrm{Rep}(Q)\mid \Hom(X,X_j)=\Ext(X,X_j)=0\text{ for }j=1,\ldots,r\},\] 
\[S^{\perp}=\{X\in\mathrm{Rep}(Q)\mid \Hom(X_j,X)=\Ext(X_j,X)=0\text{ for }j=1,\ldots,r\}.\]
It is straightforward to check that these categories are closed under direct sums, direct summands, extensions, images, kernels and cokernels. For two roots $\beta$ and $\gamma$ we denote by $\beta\in\gamma^{\perp}$ if $\hom(\gamma,\beta)=\ext(\gamma,\beta)=0$.\\
In the following we do not always distinguish between a real root and the unique indecomposable representation of this dimension. From \cite[Theorems 2.3 and 2.4]{sc2} it follows:
\begin{satz}\label{perpcat}Let $Q$ be a quiver with $n$ vertices and $S=(\alpha_1,\ldots,\alpha_r)$ be an exceptional sequence. 
\begin{enumerate}
\item The categories $^{\perp}S$ and $S^{\perp}$ are equivalent to the categories of representations of quivers $Q(^{\perp}S)$ and $Q(S^{\perp})$ respectively such that these quivers have $n-r$ vertices and no oriented cycles.
\item There is an isometry with respect to the Euler form between the dimension vectors of  $Q(^{\perp}S)$ (resp. $Q(S^{\perp})$) and the dimension vectors of $^{\perp}S$ (resp. $S^{\perp}$) given by
$\Phi((d_1,\ldots,d_{n-r}))=\sum_{i=1}^{n-r}d_i\beta_i$ where $\beta_1,\ldots,\beta_{n-r}$ are the dimension vectors of the simple representations of the perpendicular categories.
\end{enumerate}
\end{satz}
We proceed with summarising some results of \cite{rin}. For an exceptional module $S$ and a full subcategory $\mathcal{C}$ of $\Rep(Q)$ let $\mathcal{C}/S$ be the category with the same objects as $\mathcal{C}$ and the same maps modulo those factorising through $\oplus^n_{i=1} S$ for some $n\in\Nn$. We define the following full subcategories of $\Rep(Q)$:
Let $\mathcal{M}^{-S}$ be the category of representations $X$ with $\Hom(X,S)=0$ and $\mathcal{M}_{-S}$ the category of representations $X$ with $\Hom(S,X)=0$. Moreover, we define $\mathcal{M}^S$ to be the category of representations $X$ with $\Ext(S,X)=0$ such that, moreover, there does not exist a direct summand of $X$ which can be embedded into a direct sum of copies of $S$ and, finally, let $\mathcal{M}_S$ be the category of representations $X$ with $\Ext(X,S)=0$ such that, moreover, no direct summand of $X$ is a quotient of a direct sum of copies of $S$.\\
Let $X\in\mathcal{M}^S$ and $\mathcal{B}:=\{\varphi_1,\ldots,\varphi_n\}$ be a basis of $\Hom(X,S)$. Following \cite{rin} there exists an exact sequence 
\[\ses{X^{-S}}{X}{\bigoplus_{i=1}^n S}\]
induced by the basis $\mathcal{B}$ such that the induced sequences $e_1,\ldots,e_n$ form a basis of $\Ext(S,X^{-S})$. Moreover, we have $X^{-S}\in\mathcal{M}^{-S}$. Note that, equivalently, we get the representation $X^{-S}$ by taking the intersection of the kernels of all maps $X\rightarrow S$.\\
The other way around, if $Y\in\mathcal{M}^{-S}$ and $\{e_1,\ldots,e_n\}$ is a basis of $\Ext(S,Y)$ we have an induced sequence
\[\ses{Y}{Y^S}{\bigoplus_{i=1}^nS}\]
such that $Y^{S}\in\mathcal{M}^S$.\\ 
Then we have the following theorem:
\begin{satz}\label{ringel}
\begin{enumerate} 
\item  There exists an equivalence of categories given by the functor $F:\mathcal{M}^S/S\rightarrow\mathcal{M}^{-S},\,X\mapsto X^{-S}.$
\item There exists an equivalence of categories given by the functor $G:\mathcal{M}_S/S\rightarrow\mathcal{M}_{-S},\,X\mapsto X_{-S},$ where $X_{-S}=X/X'$ and $X'$ is the sum of all images of all maps $S\rightarrow X$.
\item There exist equivalences $\Psi:\mathcal{M}^{S}_{-S}\rightarrow\mathcal{M}^{-S}_{S}$ and $\Phi:\mathcal{M}^S_S/S\rightarrow\mathcal{M}^{-S}_{-S}$ induced by the above ones.
\end{enumerate}
\end{satz}
If $\alpha=\underline{\dim} S$ we define $\mathcal{M}_{\alpha}:=\mathcal{M}_S/S$. We proceed in the same manner in the other three cases.\\From \cite[Theorem 4.1]{sch} we get:
\begin{satz}\label{schofield}
Let $\alpha$ and $\beta$ be Schur roots such that $\mathrm{ext}(\alpha,\beta)=0$. Then we either have $\mathrm{hom}(\beta,\alpha)=0$ or $\mathrm{ext}(\beta,\alpha)=0$. If $\alpha$ and $\beta$ are imaginary, then we have $\mathrm{hom}(\beta,\alpha)=0$.
\end{satz}
Moreover, from \cite[Lemma 4.1]{hr} we get:
\begin{lem}\label{happelringel}
Let $X$ and $Y$ be two indecomposable representations of $Q$ such that we have $\Ext(Y,X)=0$. Then every non-zero homomorphism $f:X\rightarrow Y$ is either injective or surjective. In particular, every exceptional representation is Schurian.
\end{lem}

\subsection{On the construction of tree modules}
In this subsection we state some very useful results which can be used to construct indecomposable tree modules. Roughly speaking, given two suitable indecomposable tree modules, we consider certain exact sequences between these modules (resp. direct sums of these modules) in order to glue them appropriately. In this way, we obtain indecomposable tree modules of greater dimension vectors.
\begin{lem}\label{End}
\begin{enumerate}
\item Let $M$ and $N$ be two representations of a quiver $Q$ such that we have $\Hom(M,N)=\Hom(N,M)=0$ and $\End(N)=k$. Let $\dim_k\Ext(N,M)=d$. Let $e_1,\ldots,e_l\in\Ext(N,M)$ with $1\leq l\leq d$ be linear independent. Consider the exact sequence \[e:\ses{M}{X}{N^l}\] induced by $e_1,\ldots,e_l$. Then we have $\End(X)\subseteq\End(M)$.
\item Let $M$ and $N$ be two representations of a quiver $Q$ such that $\Hom(M,N)=\Hom(N,M)=0$ and $\End(N)=k$. Let $\dim_k\Ext(M,N)=d$. Let $e_1,\ldots,e_l\in\Ext(M,N)$ with $1\leq l\leq d$ be linear independent. Consider the exact sequence \[e:\ses{N^l}{X}{M}\] induced by $e_1,\ldots,e_l$. Then we have $\End(X)\subseteq\End(M)$.
\end{enumerate}
\end{lem}
{\it Proof.} Consider the following long exact sequence\\
\begin{xy}\xymatrix@R15pt@C20pt{0\ar[r]&\Hom(N,M)=0\ar[r]&\Hom(N,X)\ar[r]&\Hom(N,N^l)\ar[r]^{\phi}&\Ext(N,M)^1}\end{xy}\\
induced by $e$. By construction $\phi$ is injective and, therefore, $\Hom(N,X)=0$. Now consider the following commutative diagram induced by $e$:\\
\begin{xy}\xymatrix@R15pt@C20pt{&0\ar[d]&0\ar[d]&0\ar[d]\\0\ar[r]&\Hom(N^l,M)=0\ar[d]\ar[r]&\Hom(N^l,X)\ar[r]\ar[d]&\Hom(N^l,N^l)\ar[d]\\
0\ar[r]&\Hom(X,M)\ar[r]\ar[d]&\Hom(X,X)\ar[r]\ar[d]^{\phi_1}&\Hom(X,N^l)\ar[d]\\0\ar[r]&\Hom(M,M)\ar[r]^{\phi_2}&\Hom(M,X)\ar[r]&\Hom(M,N^l)=0}\end{xy}\\\\
We also have $\Hom(N^l,X)=0$. Thus $\phi_1$ is also injective and since $\phi_2$ is an isomorphism, the claim follows.\\
The dual lemma follows analogously by applying the functors $\Hom(\rule{0.2cm}{0.4pt},N)$ and $\Hom(\rule{0.2cm}{0.4pt},N^l)$ respectively.
\qed
It is easy to verify the following lemma, see also \cite[Lemma IV.1.12]{ass}:
\begin{lem}
Let
\[
\begin{xy}
\xymatrix{0\ar[r]&M\ar[r]^{\gamma}\ar[d]^{\psi}&X\ar[d]^{\phi}\ar[r]^{\delta}&N\ar[d]^{\pi}\ar[r]&0\\
0\ar[r]&M\ar[r]^{\gamma}&X\ar[r]^{\delta}&N\ar[r]&0
}\end{xy}
\]
be commutative diagram such that the rows are exact and do not split. Then we have:
\begin{enumerate}
\item If $M$ is indecomposable and $\pi$ an automorphism, then $\phi$ and $\psi$ are automorphisms.
\item If $N$ is indecomposable and $\psi$ an automorphism, then $\phi$ and $\pi$ are automorphisms.
\end{enumerate}
\end{lem}
From this we obtain the following lemma:

\begin{lem}\label{indecomp}
Let $M$ and $N$ be two indecomposable representations of a quiver $Q$ such that $\Hom(M,N)=\Hom(N,M)=0$ and $\Ext(N,M)\neq 0$. Then the representation $X$ given by some non-splitting exact sequence \[e:\ses{M}{X}{N}\] is indecomposable.
\end{lem}
{\it Proof.} Let $\phi\in\End(X)$. Since $\Hom(M,N)=0$, by the universal property of the kernel and cokernel respectively, this induces two unique endomorphisms $\psi\in\End(M)$ and $\pi\in\End(N)$. In particular, we get the following commutative diagram:
\[
\begin{xy}
\xymatrix{0\ar[r]&M\ar[r]^{\gamma}\ar[d]^{\psi}&X\ar[d]^{\phi}\ar[r]^{\delta}&N\ar[d]^{\pi}\ar[r]&0\\
0\ar[r]&M\ar[r]^{\gamma}&X\ar[r]^{\delta}&N\ar[r]&0
}\end{xy}
\]
If $\psi=0=\pi$, by the Snake lemma this induces a morphism $h:N\rightarrow M$. But since $\Hom(N,M)=0$, for instance by the preceding Lemma we get that $\mathrm{coker}(\phi)=\ker(\phi)=X$ and thus $\phi=0$. Thus we get an embedding of rings $F:\End(X)\hookrightarrow\End(M)\times\End(N)$.\\
Let $\phi\in\End(X)$ be nilpotent. Since $F(\phi)=F(\phi^2)=F(\phi)\circ F(\phi)$, we get that the induced morphisms $\pi$ and $\psi$ are nilpotent. Since $M$ and $N$ are indecomposable we have $\psi\in\{0,\mathrm{id}_M\}$ and $\pi\in\{0,\mathrm{id}_N\}$. But, because of the preceding Lemma, we get $(\psi,\pi)\in\{(\mathrm{id}_M,\mathrm{id}_N),(0,0)\}$. Thus the only nilpotent endomorphism of $X$ are $\mathrm{id}_X$ and $0_X$. Thus $\End(X)$ is local and, therefore, $X$ is indecomposable.
\qed
Let $X$ and $Y$ be two representations of a quiver $Q$. Then we can consider the linear map
\[\gamma_{X,Y}:\bigoplus_{i\in Q_0}\Hom_k(X_i,Y_i)\rightarrow\bigoplus_{\rho:i\rightarrow j\in Q_1}\Hom_k(X_i,Y_j)\]
defined by $\gamma_{X,Y}((f_i)_{i\in Q_0})=(Y_{\rho}f_i-f_jX_{\rho})_{\rho:i\rightarrow j\in Q_1}$.\\
We have $\ker(\gamma_{X,Y})=\Hom(X,Y)$ and $\mathrm{coker}(\gamma_{X,Y})=\Ext(X,Y)$, see \cite{rin2}. Whence the first assertion is easy to see, the second one follows because every exact sequence $E(f)\in\Ext(X,Y)$ is defined by a morphism $f\in\bigoplus_{\rho:i\rightarrow j\in Q_1}\Hom_k(X_i,Y_j)$ in the following way
\[\ses{Y}{((Y_i\oplus X_i)_{i\in Q_0},(\begin{pmatrix}Y_{\rho}&f_{\rho}\\0&X_{\rho}\end{pmatrix})_{\rho\in Q_1})}{X}\]
with the canonical inclusion on the left hand side and the canonical projection on the right hand side. Now it is straightforward to check that two sequences $E(f)$ and $E(g)$ are equivalent if and only if $f-g\in\mathrm{Im}(\gamma_{X,Y})$.\\\\
Let $M_{m,n}(k)$ be the set of $m\times n$ matrices with coefficients in $k$ and for $M\in M_{m,n}(k)$ denote by $M_{i,j}$ the $(i,j)$-entry. We denote by $E(s,t)\in M_{m,n}(k)$ the matrix with $E_{s,t}=1$ and zero otherwise. We call a basis $\{E(f_1),\ldots,E(f_n)\}$ of $\Ext(X,Y)$ tree-shaped if for all $i=1,\ldots, n$ we have $(f_i)_{\rho}=E(s,t)$ for exactly one $\rho\in Q_1$ and $(f_i)_{\rho '}=0$ for $\rho '\neq \rho$.\\
Since we can clearly choose a tree-shaped basis $\mathcal{B}$ of $\bigoplus_{\rho:i\rightarrow j\in Q_1}\Hom_k(X_i,Y_j)$, we can choose a basis of $\Ext(X,Y)$ consisting of elements of the form $b+\mathrm{Im}(\gamma_{X,Y})$ with $b\in\mathcal{B}$. In summary we get the following lemma: 
\begin{lem}\label{treeshapedbasis}
For every two representations $X$ and $Y$ of a quiver $Q$ there exists a tree-shaped basis of $\Ext(X,Y)$.
\end{lem}
For a real root $\alpha$ we denote the unique indecomposable representation by $X_{\alpha}$.
Now let $(X_{\alpha},X_{\beta})$ be an exceptional sequence such that $\Hom(X_{\beta},X_{\alpha})=0$. We say that a representation $Z$ has a filtration with factors $X_{\alpha}$ and $X_{\beta}$ if there exists an exact sequence
\[\ses{X_{\alpha}^e}{Z}{X_{\beta}^d}\]
with $d,e\in\Nn$. All such objects form a full subcategory $\mathcal{F}(X_{\alpha},X_{\beta})$ of $\Rep(Q)$. Moreover, it is well-known that $\mathcal{F}(X_{\alpha},X_{\beta})$ is equivalent to the category of representations of the generalised $m$-Kronecker quiver $K(m)$ with $K(m)_0=\{q_0,q_1\}$ and $K(m)_1=\{\rho_i:q_0\rightarrow q_1\mid i\in\{1,\ldots,m\}\}$ where $m=\dim\Ext(X_{\beta},X_{\alpha})$, see for instance \cite{rin2}. Then we have the following Proposition:
\begin{pro}\label{kroneckertree}
For every root of $Q$ of the form $d\beta+e\alpha$ there exists an indecomposable tree module. If $(d,e)$ is a real root of $K(m)$, then $d\beta+e\alpha$ is also a real root.

\end{pro}
{\it Proof.}  By \cite[Theorem 3.9]{wei}, for every root of $K(m)$ there exists an indecomposable tree module. Thus choose an indecomposable tree module $T_{d,e}$ of dimension $(d,e)$. Moreover, by \cite{rin2} every exceptional representation is a tree module. By Lemma \ref{treeshapedbasis} we can choose a tree-shaped basis of $\Ext(X_{\beta},X_{\alpha})$. 
Consider the exact sequence
\[\ses{X_{\alpha}^e}{Z}{X_{\beta}^d}\]
induced by $T_{d,e}$ and with respect to the chosen tree-shaped basis. Obviously, $Z$ is an indecomposable representation of dimension $d\beta+e\alpha$. Now the induced coefficient quiver has 
\[e(\dim X_{\alpha}-1)+d(\dim X_{\beta}-1)+(d+e-1)=\dim Z-1\]
vertices.\\
The second part follows from $\Sc{(d,e)}{(d,e)}=d^2+e^2-med=1$.
\qed
\begin{bem}\label{iso}
\end{bem}
\begin{itemize}
\item Consider again the generalised Kronecker quivers $K(m)$. Let $r(d,e):=(e,me-d).$ If $(d,e)\neq r^l(n,kn)$ for all $l\geq 0$ and $(n,kn)\neq(1,1)$, we even get that there exists a stable tree module for every root $(d,e)$ of $K(m)$, see \cite[Theorem 3.9]{wei}. In particular, the corresponding representation $Z$ is Schurian.\\
It can be checked easily that there only exist isotropic roots for the Kronecker quiver $K(2)$, see also \cite[Section 6]{kac}. They are given by $(d,d)$. If $d=2$, one of the two indecomposable tree modules is given by
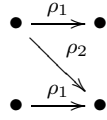
\[
\begin{xy}
\xymatrix@R20pt@C20pt{\bullet\ar[r]^{\rho_1}\ar[rd]^{\rho_2}&\bullet\\
\bullet\ar[r]^{\rho_1}&\bullet
}\end{xy}
\]
Since this is no Schur root, there does not exist any stable representation. The indecomposable tree modules of the root $(d,d)$ have the same shape. 
\end{itemize}
The next two lemmas deal with the construction of tree modules. Whence the first lemma just deals with the nature of tree-shaped bases, the second lemma states that certain submodules of indecomposable modules, which were constructed using the reflection functor, are also indecomposable modules.
\begin{lem}\label{tree}
Let $X$ be an indecomposable tree module and $S$ be an exceptional representation. Let $\mathcal{B}$ and $\mathcal{B}'$ be tree-shaped bases of $\Ext(S,X)$ and $\Ext(X,S)$ respectively. Moreover, let $\{e_1,\ldots,e_n\}\subseteq\mathcal{B}$ and $\{e'_1,\ldots,e'_k\}\subseteq\mathcal{B}'$. Consider the exact sequences induced by these bases:
\[\ses{X}{Y}{\bigoplus_{i=1}^n S}\]
and
\[\ses{\bigoplus_{i=1}^k S}{Y'}{X}.\]
Then the representations $Y$ and $Y'$ are tree modules.
\end{lem}
{\it Proof.} Obviously, the induced coefficient quiver of $Y$ has $\dim X-1+n(\dim S-1)+n=\dim Y-1$ vertices. We proceed analogously for $Y'$.
\qed
In the following we will not always state the dual lemma if it is obvious. But we will mention it if there exists one. In all these cases, the statements can be proven analogously.
\begin{lem}\label{treesub}
Let $X$ and $S$ be indecomposable such that $\Ext(S,S)=0$. % and $\Hom(X,S)=0$. 
Let $e_1,\ldots,e_n$ be a basis of $\Ext(S,X)$. Consider the exact sequence induced by this basis:
\[\ses{X}{X^S}{\bigoplus_{i=1}^n S}.\]
Moreover, consider
\[\ses{X}{Y^S}{\bigoplus_{i=1}^k S}\]
induced by $e_1,\ldots,e_k$. If $X^S$ is indecomposable, then $Y^S$ is indecomposable.\\
The dual statement of this Lemma also holds.
\end{lem}
{\it Proof.}
Consider the following commutative diagram with exact rows and columns
\[
\begin{xy}
\xymatrix@R20pt@C20pt{
&&0&0\\&0\ar[r]&\bigoplus_{i=1}^{n-k}S\ar[u]\ar@{=}[r]&\bigoplus_{i=1}^{n-k}S\ar[u]\ar[r]&0\\
0\ar[r]&X\ar[r]\ar[u]&X^S\ar[u]\ar[r]&\bigoplus_{i=1}^{n}S\ar[u]\ar[r]&0\\
0\ar[r]&X\ar@{=}[u]\ar[r]&Y^S\ar[u]\ar[r]&\bigoplus_{i=1}^{k}S\ar[u]\ar[r]&0\\
&0\ar[u]&0\ar[u]&0\ar[u]&}
\end{xy}
\]
Assume that $Y^S\cong Y_1\oplus Y_2$. Clearly, we have $\dim\Ext(S,Y^S)=n-k$ and that, by construction, $e_{k+1},\ldots,e_n$ is a basis of $\Ext(S,Y^S)$. Moreover, the sequence in the middle column is induced by this basis. Assume that $e'_{k+1},\ldots,e'_l$ is a basis of $\Ext(S,Y_1)$ and $e'_{l+1},\ldots,e'_{n}$ is a basis of $\Ext(S,Y_2)$. Then this induces an exact sequence
\[\ses{Y_1\oplus Y_2}{X_1\oplus X_2}{\bigoplus_{i=1}^{n-k}S}.\]
But then it easy to check that $X_1\oplus X_2\cong X^S$. Indeed, by construction this sequence induces an isomorphism $\Hom(\oplus_{i=1}^{n-k} S,\oplus_{i=1}^{n-k} S)\cong \Ext(\oplus_{i=1}^{n-k} S, Y_1\oplus Y_2)$. 
\qed
%Denote by $a(X)$ the minimal number of arrows of a coefficient quiver in the set of all %coefficient quivers of a representation $X$.
%\begin{lem}\label{subtree}
%Let $Z$ and $S$ be two indecomposable tree modules and let 
%\[\ses{X}{Z}{\oplus_{i=1}^t S}\]
%be a short exact sequence such that $X$ is indecomposable. Then $X$ is a tree module.\\The dual %statement of this Lemma also holds.
%\end{lem}
%{\it Proof.} Since $Z$ and $S$ are tree modules, their coefficient quivers have at most $\dim %Z-1$ and $\dim S-1$ arrows respectively. The coefficient quiver of the kernel of an epimorphism %$Z\rightarrow S$ has $\dim S$ points less, and if $Z$ is indecomposable, at least $\dim S$ %arrows less.
%Thus we have 
%\[\dim X+t\dim S-1= a(Z)\geq a(X)+t\dim S .\]
%Thus we have $a(X)=\dim X-1$ because $X$ is indecomposable.
%\qed
Let $X_{\alpha}$ be an exceptional representation and $Y\in X_{\alpha}^{\perp}$ with $\underline{\dim}Y=\beta$ such that $Y$ does not embed into a direct sum of copies of $X_{\alpha}$. Note that if $\beta$ is an imaginary root and $Y$ is indecomposable, this is automatically satisfied. Indeed, if we had $Y\hookrightarrow\oplus X_{\alpha}$, the corresponding long exact sequence would induce an epimorphism $\Ext(\oplus X_{\alpha},Y)\twoheadrightarrow\Ext(Y,Y)$. But since $\Ext(Y,Y)\neq 0$, this is not possible because $Y\in X_{\alpha}^{\perp}$.
 
\begin{pro}\label{treeconstruction}
\begin{enumerate}\item If $Y$ is an indecomposable tree module such that $Y\in X_{\alpha}^{\perp}$ and $\underline{\dim} Y=\beta$, there exists an indecomposable tree module of dimension $\beta+r\alpha$ for all $0\leq r\leq \dim\Ext(Y,X_{\alpha})$. \\
Moreover, if $Y$ is Schurian and $\Hom(Y,X_{\alpha})=0$, then there exists a Schurian tree module of dimension $\beta+r\alpha$ for all $0\leq r\leq \dim\Ext(Y,X_{\alpha})$.\\
The dual statement also holds.
\item If $Y$ is an indecomposable tree module such that $Y\in\mathcal{M}^{-\alpha}_{-\alpha}$ with $\underline{\dim} Y=\beta$, there exists an indecomposable tree module of dimension $\beta+r\alpha$ for all $0\leq r\leq \dim\Ext(Y,X_{\alpha})+\dim\Ext(X_{\alpha},Y)$. If $Y$ is Schurian, there exists an Schurian tree module of dimension $\beta+r\alpha$ for all $0\leq r\leq \max\{\dim\Ext(Y,X_{\alpha}),\dim\Ext(X_{\alpha},Y)\}.$
\end{enumerate}
%\item If $Y$ is an indecomposable tree module such that $Y\in ^{\perp}X_{\alpha}$ and $\underline{\dim} Y=\beta$, there exists an indecomposable tree module of dimension $\beta+r\alpha$ for all $-\dim\Hom(X_{\alpha},Y)\leq r\leq \dim\Ext(X_{\alpha},Y)+\dim\Hom(X_{\alpha},Y)$. \\
%Moreover, if $Y$ is Schurian, then there exists a Schurian tree module of dimension %$\beta+r\alpha$ for all $-\dim\Hom(X_{\alpha},Y)\leq r\leq\max\{0,-\Sc{\alpha}{\beta}\}$.
%\end{enumerate}
\end{pro}
{\it Proof.}
We first prove the first statement. If $\beta$ is a real Schur root, we either have $\Hom(Y,X_{\alpha})=0$ or $\Ext(Y,X_{\alpha})=0$. Indeed, a general representation of dimension $\beta$ is Schurian, and thus we can apply Theorem \ref{schofield}. In the first case, we are in the situation of Proposition \ref{kroneckertree} and Lemma \ref{End}. In the second case, the statement is trivial.\\\\
Thus let $\beta$ be no real Schur root. Since $Y$ does not embed into a direct sum of copies of $X_{\alpha}$, we have $Y\in\mathcal{M}^{\alpha}_{-\alpha}$. Thus by Theorem \ref{ringel} and Lemma \ref{tree} we have that $Y^{\alpha}$, which is given by the short exact sequence
\[\ses{\bigoplus_{i=1}^{\dim\Ext(Y,X_{\alpha})}X_{\alpha}}{Y^{\alpha}}{Y}\]
induced by a tree-shaped basis of $\Ext(Y,X_{\alpha})$, is an indecomposable tree-module. Now, by applying Lemmas \ref{tree} and \ref{treesub}, we obtain indecomposable factor modules of $Y^{\alpha}$ which are trees for all $0\leq r\leq \dim\Ext(Y,X_{\alpha})$. For instance by Lemma \ref{End}, we obtain that if $\Hom(Y,X_{\alpha})=0$ and $Y$ is Schurian, the constructed representations are Schurian.\\\\
%If $Y$ is Schurian and $\Sc{\beta}{\alpha}\leq 0$, we choose a tree-shaped basis of $\Ext(Y^{-\alpha},X_{\alpha})\cong\Ext(Y,X_{\alpha})$. Then we can apply Lemma \ref{End} in order to get a Schurian tree module of dimension $\beta+r\alpha$ for all  $-\dim\Hom(Y,X_{\alpha})\leq r\leq -\Sc{\beta}{\alpha}$. If $\Sc{\beta}{\alpha}\geq 0$, we choose a tree-shaped basis of $\Ext(X_{\alpha},Y^{-\alpha})$ and proceed analogously to the other case.\\\\
The second statement is obtained in the following way: let $n=\dim\Ext(Y,X_{\alpha})$ and $m=\dim\Ext(X_{\alpha},Y)$. Consider the commutative diagram with exact rows and columns induced by tree-shaped bases of $\Ext(Y,X_{\alpha})$ and $\Ext(X_{\alpha},Y)$ respectively:
\[
\begin{xy}
\xymatrix@R20pt@C20pt{
&0\ar[d]&0\ar[d]&0\ar[d]\\0\ar[r]&\bigoplus_{i=1}^{n}X_{\alpha}\ar[d]\ar@{=}[r]&\bigoplus_{i=1}^{n}X_{\alpha}\ar[d]\ar[r]&0\ar[r]\ar[d]&0\\
0\ar[r]&Y^{\alpha}\ar[r]\ar[d]&Y_{\alpha}^{\alpha}\ar[d]\ar[r]&\bigoplus_{i=1}^{m}X_{\alpha}\ar@{=}[d]\ar[r]&0\\
0\ar[r]&Y\ar[d]\ar[r]&Y_{\alpha}\ar[d]\ar[r]&\bigoplus_{i=1}^{m}X_{\alpha}\ar[d]\ar[r]&0\\
&0&0&0&}
\end{xy}
\]
The equivalence $\Phi^{-1}:\mathcal{M}^{-\alpha}_{-\alpha}\rightarrow\mathcal{M}^{\alpha}_{\alpha}$ is given by $\Phi^{-1}(Y)=Y^{\alpha}_{\alpha}$. Thus by Theorem \ref{ringel} and Lemma \ref{tree}, we have that $Y_{\alpha}^{\alpha}$ is an indecomposable tree module of dimension $\beta+(n+m)\alpha$.\\
Now by applying Lemmas \ref{tree} and \ref{treesub}, we can also construct indecomposable subtrees and factor trees of $Y_{\alpha}^{\alpha}$, $Y^{\alpha}$ and $Y_{\alpha}$ respectively such that the claim follows.
\qed
\begin{bem}
\end{bem}
\begin{itemize}
\item We have $\dim\End(Y_{\alpha}^{\alpha})=\dim\End(Y)+\Sc{\underline{\dim}Y^{-\alpha}}{\underline{\dim}X_{\alpha}}\Sc{\underline{\dim}X_{\alpha}}{\underline{\dim}Y^{-\alpha}}$ by \cite[Proposition $3^{\ast}$]{rin}. Using the same arguments, we can calculate the dimension of the endomorphism rings of all other tree modules considered in the preceding Proposition.
\item The preceding Proposition also gives a recipe how to construct indecomposable tree modules for non-Schur roots. In particular, given a Schur root $\beta$ the idea is to construct a tree module such that $X_{\beta}\in\alpha^{\perp}$ (resp. $X_{\beta}\in$$^{\perp}\alpha$) for a real Schur root $\alpha$ and, moreover, $\Hom(X_{\beta},X_{\alpha})\neq 0$ and $\Ext(X_{\beta},X_{\alpha})\neq 0$ (resp. $\Hom(X_{\alpha},X_{\beta})\neq 0$ and $\Ext(X_{\alpha},X_{\beta})\neq 0$). See also Section \ref{examples} for a more detailed discussion of an example.
\end{itemize}
The following Proposition is based on the algorithm of \cite[Section 4]{dw}. Roughly speaking, it describes what kind of possibilities there are to decompose an imaginary Schur root into smaller Schur roots. Later on, this decomposition will be one of the basic tools for the construction of tree modules.
\begin{pro}\label{zerlegung}
Let $\alpha$ be an imaginary Schur root. Then at least one the following cases holds:
\begin{enumerate}
\item There exist a real Schur root $\beta$ and $t\in\Nn_+$ such that $\gamma=\alpha-t\beta$ is an imaginary Schur root. Moreover, we have $\beta\in\gamma^{\perp}$ and $\hom(\beta,\gamma)=0$ or $\beta\in$$^{\perp}\gamma$ and $\hom(\gamma,\beta)=0$.
\item There exist a real Schur root $\beta$ and a real or isotropic root $\gamma$ and $d,e\in\Nn_+$ such that $\alpha=\beta^d+\gamma^e$. Moreover,  we have $\beta\in\gamma^{\perp}$ and $\hom(\beta,\gamma)=0$ or $\beta\in$$^{\perp}\gamma$ and $\hom(\gamma,\beta)=0$ and $(d,e)$ is a root of $K(\ext(\beta,\gamma))$ or $K(\ext(\gamma,\beta))$.
\item There exist two imaginary Schur roots $\gamma$ and $\delta$ such that $\gamma+\delta=\alpha$. Moreover, we have $\delta\in\gamma^{\perp}$ and $\hom(\delta,\gamma)=0$. 
\end{enumerate}
\end{pro}
{\it Proof.} The statement follows when considering the algorithm of Derksen and Weyman \cite[Section 4]{dw}.
\qed
From the preceding Proposition we obtain the following Corollary:
\begin{kor}\label{isotropic}
For every isotropic root $\alpha$ there exists a decomposition 
$\alpha=\beta^{kd}+\gamma^{d}$ where $\beta$ is a real Schur root such that $\gamma\in\beta^{\perp}$ or $\gamma\in$$^{\perp}\beta$. Additionally, $\gamma$ is either an isotropic or a real Schur root such that we have $\mathrm{hom}(\gamma,\beta)=0$ or $\mathrm{hom}(\beta,\gamma)=0$. If $\alpha$ is indivisible, then $\gamma$ is indivisible such that if $\gamma$ is isotropic, we have $d=1$, and if $\gamma$ is real, we have $d=1$ and $k=1$.
\end{kor}
{\it Proof.}
The canonical decomposition of $\alpha$ is given by $\alpha=\oplus\tilde{\alpha}$ for some indivisible isotropic root $\tilde{\alpha}$. We first assume that $\alpha$ is indivisible. There cannot be a decomposition $\alpha=\beta+\gamma$ such that $\gamma$ and $\beta$ are imaginary and, moreover, $\gamma\in\beta^{\perp}$, $\ext(\gamma,\beta)\neq 0$ and $\hom(\gamma,\beta)= 0$. Indeed, otherwise we had
\[\Sc{\alpha}{\alpha}=\Sc{\beta}{\beta}+\Sc{\gamma}{\gamma}+\Sc{\gamma}{\beta}<0.\]
Thus by Proposition \ref{zerlegung} we have a decomposition $\alpha=\beta^e+\gamma^d$ such that without loss of generality we have $\gamma\in\beta^{\perp}$ and, moreover, $\beta$ is real and  $\gamma$ is either imaginary or real. But since
\[\Sc{\gamma^d}{\gamma^d}=\Sc{\alpha-e\beta}{\gamma^d}=\Sc{\alpha}{\alpha}-e\Sc{\alpha}{\beta}\]
and $\Sc{\alpha}{\beta}\leq 0$, we get that $\gamma$ is isotropic or real. Indeed, since a general representation of dimension $\alpha$ is Schurian and has a subrepresentation of dimension $\beta$ we have $\mathrm{hom}(\alpha,\beta)=0$.\\
If $\gamma$ is isotropic, we get $\Sc{\alpha}{\beta}=0$. Thus $\gamma=\alpha-\Sc{\beta}{\alpha}\beta$ is a Schur root and we obtain $d=1$ and $k=\Sc{\gamma}{\beta}=-\Sc{\beta}{\alpha}$. In particular, $\gamma$ is indivisible because otherwise $\alpha=s\tilde{\gamma}+s\Sc{\tilde{\gamma}}{\beta}\beta$ would be divisible as well.\\
If $\gamma$ is real, it follows that $\Sc{\alpha}{\beta}=-d$, $e=d$ and $k=1$. Indeed, there only exist isotropic roots for $K(m)$ if $m=2$ and they are given by $(d,d)$, see Remark \ref{iso}. In particular, since, by assumption, $\alpha$ is indivisible, we obtain $d=1$.\\
If $\alpha$ is divisible, say $\alpha=d\tilde{\alpha}$, it is now easy to check that we get a decomposition $\alpha=\beta^{dk}+\gamma^d$. Thereby, as before we get a decomposition $\tilde{\alpha}=\beta^k+\gamma$.
\qed
\subsection{Trees of isotropic roots}\label{secisotropic}
In this section we construct indecomposable tree modules for every isotropic root.
By Corollary \ref{isotropic} we have a decomposition $d\alpha=\beta^{dk}+\gamma^d$ where $\alpha$ is an indivisible isotropic root, $\beta$ is a real Schur root and $\gamma$ is either an isotropic or a real Schur root. Because of Theorem \ref{perpcat} after at most $|Q_0|-1$ steps we get that $\gamma$ is real. Indeed, $\gamma$ corresponds to a root of a quiver with one vertex less, and for the quiver with only one vertex and without arrows there exists only one indecomposable representation. Thus let us first assume that $\gamma$ is real. Then by Corollary \ref{isotropic} we get $k=1$ and $\dim\Ext(X_{\gamma},X_{\beta})=2$. Since $X_{\gamma}\in X_{\beta}^{\perp}$ and $\Hom(X_{\gamma},X_{\beta})=0$, the category of short exact sequences of the form 
\[\ses{X_{\beta}^d}{X_{d\alpha}}{X_\gamma^d}\]
corresponds to the representations of the isotropic root $(d,d)$ of $K(2)$. But since there exists an indecomposable tree module for $(d,d)$, by Proposition \ref{kroneckertree} there exists an indecomposable tree module $X_{d\alpha}$.\\
Thus, in general, we get a sequence $(\alpha_n,\ldots,\alpha_1)$ of real Schur roots such that $\sum_{i=1}^n\alpha_i^{k_i}=\alpha$. Applying successively Proposition \ref{treeconstruction} and starting as mentioned above, we can inductively construct indecomposable tree modules $Y_l$ with $Y_l\in\alpha_i^{\perp}$ or $Y_l\in$$^{\perp}\alpha_i$ and $\underline{\dim} Y_l=\sum_{j=1}^l\alpha_j^{k_j}$ for $n\geq i>l$ and $l=1,\ldots,n$. In summary we obtain the following result:
\begin{satz}\label{isot}
For every isotropic root of a quiver $Q$ there exist more than one isomorphism class of indecomposable tree modules. 
\end{satz}

\subsection{Tree modules of imaginary Schur roots}
In this section we construct indecomposable tree modules for every imaginary Schur root.
\begin{satz}\label{imag}
Let $\alpha$ be a Schur root. Then there exists an indecomposable representation $X_{\alpha}$ which is a tree module. If $\alpha$ is an imaginary Schur root, there exists more than one isomorphism class. 
\end{satz}
{\it Proof.} If $\alpha$ is a real Schur root, the claim follows from \cite{rin1}. Thus let $\alpha$ be an imaginary Schur root. We consider Proposition \ref{zerlegung}. Assume that there exists a real Schur root $\beta$ and a decomposition $\alpha=\beta^d+\gamma$ such that $\gamma$ is an imaginary Schur root or $\alpha=\beta^d+\gamma^e$ where $\gamma$ is a real Schur root or an isotropic root. Without loss of generality, for this decomposition we have 
\[\mathrm{ext}(\beta,\gamma)=\mathrm{hom}(\beta,\gamma)=\mathrm{hom}(\gamma,\beta)=0.\]
If $\gamma$ is real, $(d,e)$ is an imaginary Schur root of $K(\mathrm{ext}(\gamma,\beta))$ and we can apply Proposition \ref{kroneckertree}.\\
Now assume that $\gamma$ is imaginary. If $\gamma$ is an isotropic root, we proceed as mentioned in the preceding subsection in order to construct a tree module $X_{\gamma^e}$ of the root $\gamma^e$. More detailed, by Theorem \ref{perpcat} we have that the category $X_{\beta}^{\perp}$ is equivalent to the category of representations of a quiver with $|Q_0|-1$ vertices. Moreover, $\gamma$ corresponds to an isotropic root of this quiver. Thus by induction hypothesis there exists an indecomposable tree module $X_{\gamma}\in X_{\beta}^{\perp}$ of dimension $\gamma$. Note that we decompose $\gamma$ in $X_{\beta}^{\perp}$, but, construct the tree module using the corresponding dimension vectors and tree modules of the original quiver. Then we can apply Proposition \ref{treeconstruction} since $X_{\gamma^e}\in X_{\beta}^{\perp}$. \\If $\gamma$ is imaginary, we proceed as in the case before. Thus applying Proposition \ref{treeconstruction}, we get an indecomposable tree module of dimension $\alpha$. Moreover, it is easy to check that in both cases the construction is not unique. Thus there exists more than one isomorphism class of tree modules. \\
Note that, if $|Q_0|=2$ we deal with the $m$-Kronecker quiver $K(m)$. In this case we always get a decomposition of $\alpha$ into two real Schur roots.\\
Thus it remains to consider the last case of Proposition \ref{zerlegung}. Therefore, let $\alpha=\beta+\gamma$ be a decomposition into imaginary Schur roots such that $\gamma\in\beta^{\perp}$. By induction hypothesis there exist indecomposable tree modules $X_{\beta}$ and $X_{\gamma}$ of dimension $\beta$ and $\gamma$ such that $X_{\gamma}\in X_{\beta}^{\perp}$. Note that every Schur sequence can be refined to an exceptional sequence, see \cite[Theorem 4.11]{dw2}. Thus by Theorem \ref{perpcat}, we have that $\beta$ and $\gamma$ respectively correspond to imaginary Schur roots of two quivers with less than $|Q_0|$ vertices. Moreover, we have $\Hom(X_{\gamma},X_{\beta})=0$. Indeed, by Lemma \ref{happelringel} every homomorphism is injective or surjective. Therefore, we either get an epimorphism $\Ext(X_{\beta},X_{\gamma})\rightarrow\Ext(X_{\gamma},X_{\gamma})$ or an epimorphism $\Ext(X_{\beta},X_{\gamma})\rightarrow\Ext(X_{\beta},X_{\beta})$ and thus a contradiction. Consider some non-splitting exact sequence induced by a basis element of a tree-shaped basis, i.e.
\[\ses{X_{\beta}}{X_{\gamma+\beta}}{X_{\gamma}}.\]
Now by Lemma \ref{indecomp}, we have that $X_{\gamma+\beta}$ is indecomposable.

\qed

\section{Examples}\label{examples}

\begin{bei}\label{bei}
\end{bei}
Proposition \ref{treeconstruction} and Theorem \ref{ringel} respectively give a recipe for the construction of tree modules for every root based on the reflection functors, in particular for non-Schur roots. But this does not always work, see \cite{wie} for a counterexample.\\ We will state two examples of the subspace quiver, a negative and a positive one and proceed analogously to \cite{wie}.\\Given a non-Schur root $\alpha$ with the canonical decomposition $\alpha=\oplus\gamma_i^{k_i}$ we know that a general representation of dimension vector $\alpha$ has a factor and a subrepresentation of dimension $\gamma_i$ for all $i$. In particular, we have $\Sc{\alpha}{\gamma_i}>0$ and $\Sc{\gamma_i}{\alpha}>0$ for all $i$. Thus the real roots under these roots are possible candidates for the application of Proposition \ref{treeconstruction}.\\
Consider the $8$-subspace quiver, i.e. $Q_0=\{i_0,i_1,\ldots,i_8\}$ and $Q_1=\{\rho_j:i_j\rightarrow i_0\mid j\in\{1,\ldots,8\}\}$ and the following real root with its canonical decomposition:
\begin{eqnarray*}\alpha&=&(48,1,1,1,15,15,18,18,46)\\&=&(39,1,1,1,12,12,15,15,37)\oplus(3,0,0,0,1,1,1,1,3)^3=\alpha_1\oplus\alpha_2^3.\end{eqnarray*}
In order to find a real Schur root $\beta$ such that we have $X_{\alpha}\in X_{\beta}^{\beta}$ we
have the necessary conditions $\beta<\alpha$ and $\Sc{\beta}{\alpha}\geq 0$ and $\Sc{\alpha}{\beta}\geq 0$. Following \cite[Section 6]{sch}, see also \cite{wie} for a more detailed discussion of an example, this implies $s_{\alpha}(\beta)=s_{i_1}\ldots s_{i_n}(\beta)<0$ where $s_{\alpha}:\Zn Q_0\rightarrow\Zn Q_0$ is the reflection at the hyperplane perpendicular to $\alpha$. But this is equivalent to 
\[\beta=s_{i_n}\ldots s_{i_{k+1}}(e_{i_k})\]
for some $k\geq 1$. Now it is straightforward to determine all such roots. It turns out that the only candidate for a reflection is $\alpha_2$.\\
Thus assume that we have $X_{\alpha}\in X^{\alpha_2}_{\alpha_2}$ for the unique indecomposable representation $X_{\alpha}$. Then for the unique indecomposable representation $X_{\delta}$ with $\delta=(3,1,1,1,0,0,3,3,1)$ we would have $\Hom(X_{\alpha_2},X_{\delta})=0$. But obviously the indecomposable representation of dimension vector $(1,0,0,0,0,0,1,1,1)$ is a factor of $X_{\alpha_2}$ and a subrepresentation of $X_{\delta}$. This is a contradiction.\\
Note that this also shows that it is not possible to construct tree modules of every root via the reflection functor when restricting to the universal cover of some quiver. Indeed, for the subspace quiver every component of the universal cover looks like the subspace quiver itself. This holds for every tree-shaped quiver.
\\\\
It is easy to construct tree modules of non-Schur roots. For instance we can start with $\alpha=(1,0,0,1,1,1)$ and $\beta=(1,1,1,0,0,0)$ of the $5$-subspace quiver and the corresponding unique indecomposable representations which are obviously tree modules. We have $\dim\Ext(X_{\alpha},X_{\beta})=2$ and
$\dim\Ext(X_{\beta},X_{\alpha})=1$. Thus choosing two tree-shaped bases of the groups of extensions, we get tree modules of dimension vectors $(4,1,1,3,3,3)$ and $(4,3,3,1,1,1)$ respectively. 
\begin{bei}
\end{bei}
We consider the quiver
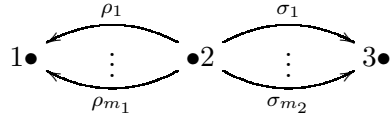
\[
\begin{xy}
\xymatrix@R10pt@C20pt{
\\1\bullet &\vdots&\bullet 2\ar@/_1.0pc/[ll]_{\rho_{1}}\ar@/^1.0pc/[ll]^{\rho_{m_1}}\ar@/_1.0pc/[rr]_{\sigma_{m_2}}\ar@/^1.0pc/[rr]^{\sigma_1}&\vdots&3\bullet\\
}
\end{xy}\]
Let $(e,d,f)$ be a Schur root of $Q$. Following the algorithm of \cite{dw} we first consider the canonical decomposition of $(e,d)$ corresponding to the generalised Kronecker quiver $K(m_1)$ with $m_1$ arrows. If $(e,d)$ is an imaginary root or a real root, we would have found a decomposition into an imaginary root or a real root and a multiple of a real root which would be the simple root in this case. Thus assume that $(e,d)$ is no root. In this case, the canonical decomposition looks like
\[(e,d)=(e_1,d_1)^{r_1}\oplus(e_2,d_2)^{r_2}.\] 
where $(e_1,d_1)$ and $(e_2,d_2)$ are real roots. Now we consider the canonical decomposition of the dimension vector $(r_1,f)$ of the generalised Kronecker quiver with $m_1d_1$ arrows. Indeed, we have $\ext((e_1,d_1,0),(0,0,1))=m_1d_1$. Again, if this is a Schur root we get a decomposition into Schur roots as claimed in Proposition \ref{zerlegung}. Otherwise we again get a decomposition into multiples of real Schur roots. If we proceed like this, after finitely many steps we get a decomposition into either multiples of two real roots or an imaginary root and a multiple of a real root.\\\\
Let us consider the example $m_1=2$, $m_2=2$ and $(e,d,f)=(7,4,5)$ which is an imaginary Schur root. We have the decomposition $(7,4)=(3,2)\oplus (2,1)^2$. Now we have to consider the dimension vector $(1,5)$ of $K(4)$. We obtain $(1,5)=(1,4)\oplus (0,1)$. Next we have to consider the dimension vector $(1,2)$ of $K(2)$ which is a real root. Thus we found the decomposition 
\[(e,d,f)=(3,2,4)+(4,2,1)\]
with $\ext((4,2,1),(3,2,4))=8$. For instance, we get the following coefficient quiver:
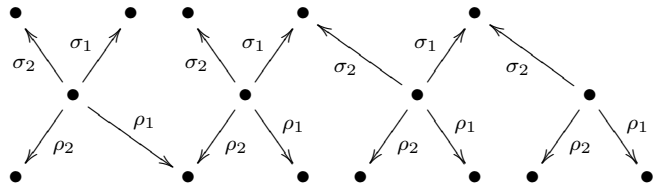
\[
\begin{xy}
\xymatrix@R20pt@C10pt{\bullet&&\bullet&\bullet&&\bullet&&&\bullet&&&\\
&\bullet\ar[lu]^{\sigma_2}\ar[ru]^{\sigma_1}\ar[ld]^{\rho_2}\ar[rrd]^{\rho_1}&&&\bullet\ar[lu]^{\sigma_2}\ar[ru]^{\sigma_1}\ar[rd]^{\rho_1}\ar[ld]^{\rho_2}&&&\bullet\ar[llu]^{\sigma_2}\ar[ru]^{\sigma_1}\ar[ld]^{\rho_2}\ar[rd]^{\rho_1}&&&\bullet\ar[llu]^{\sigma_2}\ar[rd]^{\rho_1}\ar[ld]^{\rho_2}\\
\bullet&&&\bullet&&\bullet&\bullet&&\bullet&\bullet&&\bullet&&
}
\end{xy}\]
Choosing another basis element of a tree-shaped basis we obtain the coefficient quiver
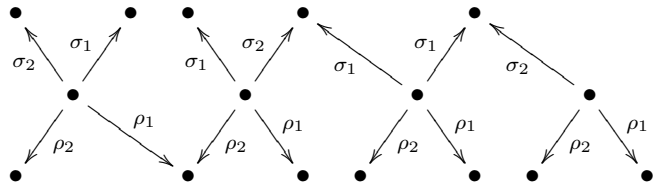
\[
\begin{xy}
\xymatrix@R20pt@C10pt{\bullet&&\bullet&\bullet&&\bullet&&&\bullet&&&\\
&\bullet\ar[lu]^{\sigma_2}\ar[ru]^{\sigma_1}\ar[ld]^{\rho_2}\ar[rrd]^{\rho_1}&&&\bullet\ar[lu]^{\sigma_1}\ar[ru]^{\sigma_2}\ar[rd]^{\rho_1}\ar[ld]^{\rho_2}&&&\bullet\ar[llu]^{\sigma_1}\ar[ru]^{\sigma_1}\ar[ld]^{\rho_2}\ar[rd]^{\rho_1}&&&\bullet\ar[llu]^{\sigma_2}\ar[rd]^{\rho_1}\ar[ld]^{\rho_2}\\
\bullet&&&\bullet&&\bullet&\bullet&&\bullet&\bullet&&\bullet&&
}
\end{xy}\]
The first representation is obviously a representation of the universal cover such that the entries of the dimension vector are either one or zero. The second one can also be seen as a representation of the universal cover by identifying vertices inducing the same word, e.g. the two terminating vertices of the two arrows coloured by $\sigma_1$ starting at the same vertex, see also Section \ref{allg}.\\\\
{\bf Acknowledgements:} I would like to thank Klaus Bongartz, Markus Reineke and Claus Michael Ringel for helpful discussions and useful comments.


\begin{thebibliography}{1}
\bibitem{ass}
Assem, I., Simson, D., Skowronski, A.: Elements of the Representation Theory of Associative Algebras. Cambridge University Press, Cambridge 2007.  
\bibitem{dw} Derksen, H., Weyman, J.: On the canonical decomposition of quiver representations. Compositio Mathematica {\bf 133}, 245-265 (2002).
\bibitem{dw3} Derksen, H., Weyman, J.: Quiver representations. Notices of the AMS {\bf 52} 2, 200-206 (2005)
\bibitem{dw2} Derksen, H., Weyman, J.: Combinatorics of quiver representations. Preprint 2007, arXiv: math/0608288v2.
\bibitem{gab} Gabriel, P.: The universal cover of a finite-dimensional algebra. Representations of algebras. Lecture Notes in Mathematics {\bf 903}, 68-105 (1981).
\bibitem{gre} Green, E.L.: Group-graded algebras and the zero-relation problem. Representations of Algebras {\bf 903}, 106-115 (1981).
\bibitem{hr} Happel, D., Ringel, C.M.: Tilted Algebras. Transactions of the American Mathematical Society {\bf 274}, no.2, 399-443 (1982).
\bibitem{kac} Kac, V.: Infinite root systems, representations of graphs and invariant theory II. Journal of Algebra {\bf 78}, 141-162 (1982).
\bibitem{rin} Ringel, C.M.: Reflection functors for hereditary algebras. J. London Math. Soc. {\bf 2}, 465-479 (1980).
\bibitem{rin1} Ringel, C.M.: Exceptional modules are tree modules. Linear algebra and its Applications {\bf 275-276}, 471-493 (1998).
\bibitem{rin3}
Ringel, Claus Michael: Combinatorial Representation Theory. History and future. Representations of Algebras {\bf 1}. Proceedings of the Conference on ICRA IX, Beijing 2000, Beijing Normal University Press, 122-144 (2002).
\bibitem{rin2} Ringel, C.M.: Representations of $K$-species and bimodules. Journal of Algebra {\bf 41}, 269-302 (1976).
\bibitem{sch} Schofield, A.: General representations of quivers. Proc. London Math. Soc. (3) {\bf 65}, 46-64 (1992).
\bibitem{sc2} Schofield, A.: Semi-invariants of quivers, J. London Math. Soc. {\bf 43}, 383-395 (1991).
\bibitem{wei} Weist, T.: Tree modules of the generalised Kronecker quiver. Journal of Algebra {\bf 323}, 1107-1138 (2010).
\bibitem{wie} Wiedemann, M.: A remark on the constructibility of real root representations of quivers using universal extension functors. Journal of Algebra {\bf 321}, 1711-1718 (2009).
\end{thebibliography}
\end{document}